\input amstex
\input amssym
\documentstyle{amsppt}
\magnification=1200
%\pageheight{12in} \pagewidth{7.5in}

\def\a{\alpha}

\def\G{\Gamma}
\def\g{\gamma}

\def\s{\sigma}

\def\x{\times}

\def\o{\overline}
\def\f{\flushpar}

\def\om{\omega}
\def\Om{\Omega}
\def\B{\Cal B}

\def\C{\Cal C}

\def\({\biggl(}
\def\){\biggr)}

\def\mpt{\text{measure preserving transformation}}
\def\eppt{\text{ergodic, probability preserving transformation}}

\def\eppG{\text{ergodic, probability preserving  $\Bbb G$-action}}\def\ppG{\text{probability preserving  $\Bbb G$-action}}
\def\ensg{\text{ergodic, non-singular $\Gamma$-action}}
\def\nsg{\text{non-singular $\Gamma$-action}}
\def\ensG{\text{ergodic, non-singular $\Bbb G$-action}}
\def\nsG{\text{non-singular $\Bbb G$-action}}
\def\<{\bold\langle}
\def\>{\bold\rangle}

\def \r{\Cal R}

\def\rwgp{\text{\tt RW}(\Bbb
G,p)} \topmatter
\title Exactness of Rokhlin endomorphisms and weak mixing of Poisson boundaries.\endtitle
\author Jon Aaronson \ \&\ Mariusz
 Lema\'nczyk\endauthor\address[Aaronson]{\ \ School of Math. Sciences, Tel Aviv University,
69978 Tel Aviv, Israel.}
\endaddress
\email{aaro\@tau.ac.il}\endemail\address[Lema\'nczyk]{\ \ Faculty
of Mathematics and Computer Science, Nicolaus Copernicus
University, ul. Chopina 12/18, 87-100 Toru\'n,
Poland}\endaddress\email{mlem\@mat.uni.torun.pl}\endemail
\abstract We give conditions for the exactness of Rokhlin
endomorphisms, apply these to random walks on locally compact,
second countable topological groups and obtain that the action on
the Poisson boundary of an adapted random walk on such a group is
 weakly mixing.
\endabstract\thanks{\copyright 2004. Aaronson was partially supported by
the EC FP5, IMPAN-BC Centre of Excellence when this work was
begun. Lema\'nczyk's research was partially supported by KBN grant
1 P03A 03826. Both authors thank MPIM, Bonn for hospitality when
this work was finished. }\endthanks\subjclass 37A20,\ 37A15
(37A40,
 60B15, 60J50)\endsubjclass
\endtopmatter
\rightheadtext{ exactness}
\document
\heading{ \S0 Introduction}\endheading\subheading{Rokhlin
endomorphisms}
 \par By a {\it non-singular endomorphism}  we mean a quadruple $(X,\B,m,T)$ where
 \f $(X,\B,m)$ is a standard probability space and $T:X_0\to X_0$ is a measurable transformation of
 $X_0\in\B,\ m(X_0)=1$ satisfying $m(T^{-1}A)=0\ \Leftrightarrow\ m(A)=0\ \
 (A\in\B)$. As in \cite{Ro2}, the endomorphism $T$ is called {\it exact} if $\goth
 T(T):=\bigcap_{n=0}^\infty
 T^{-n}\B\overset{m}\to{=}\{\emptyset,X\}$.
 \par  Let $\Bbb G$ be a locally
compact, Polish topological (LCP) group.

By a {\it non-singular $\Bbb G$-action} on a probability space
$(Y,\C,\nu)$ we mean a measurable homomorphism $S:\Bbb G\to
\text{\rm Aut}\,(Y)$ where $\text{\rm Aut}\,(Y)$ denotes the group
of invertible, non-singular transformations of $Y$ equipped with
its usual Polish topology. The action $S$ is called {\it
probability preserving} if each $S_g$ preserves $\nu$.

Given a non-singular endomorphism $(X,\B,m,T)$,  a $\nsG$ $S:\Bbb
G\to \text{\rm Aut}\,(Y)$ of a LCP group $\Bbb G$ and a measurable
function $f:X\to \Bbb G$, we consider the {\it Rokhlin
endomorphism} $\widetilde T=\widetilde T_{f,S}:X\x Y\to X\x Y$
defined by
$$\widetilde T(x,y):=  (Tx,S_{f(x)}y).$$
Rokhlin endomorphisms  first appeared in \cite{Ro1} (see also
\cite{AR}). We give conditions for their exactness  (theorem 2.3).

These conditions are applied to random walk endomorphisms.
Meilijson (in \cite{Me}) gave sufficient conditions for exactness
for random walk endomorphisms over $\Bbb G=\Bbb Z$. We clarify
Meilijson's theorem proving a converse (proposition 4.2),  extend
it to  LCP Abelian groups (theorem 4.1), characterize the
exactness of the Rokhlin endomorphism  for a steady random walk
 (theorem 4.5) and obtain
that the group action on the Poisson boundary (see \S4) of an
adapted (i.e. globally supported) random walk is weakly mixing
(proposition 4.4).

Tools employed include the ergodic theory of ``associated actions"
(see \S1), and the boundary theory of random walks (see \S4).

 The
authors would like to thank the referee for some useful comments
including a simplification in the proof of proposition 3.4.
\heading{ \S1 Associated actions}\endheading
\par For a non-singular endomorphism  $(Z,\Cal D,\nu,R)$
set \f$\bullet\ \ $   $\Cal I(R):=\{A\in\Cal D:\ R^{-1}A=A\}$ --
the {\it invariant} $\s$-algebra, and \f$\bullet\ \ $   $\Cal
T(R):=\bigcap_{n=0}^\infty R^{-n}\Cal D$ -- the {\it tail}
$\s$-algebra.

\par Let $(X,\B,m,T)$ be a non-singular endomorphism. Let $\Bbb G$ be a locally
compact, Polish topological (LCP) group, let $f:X\to \Bbb G$ be
measurable.
\par There  are two {\it
associated (right) $\Bbb G$-actions} arising from the invariant
and tail $\s$-algebras of $T_f$, which are defined as follows:
\f$\bullet\ \ $ define the (left) {\it skew product} endomorphism
$T_f:X\x\Bbb G\to X\x\Bbb G$ by $T_f(x,g):=(T x,f(x)g)$ and fix
$\Bbb P\in\Cal P(X\x\Bbb G),\ \Bbb P\sim m\x m_{\Bbb G}$;
\f$\bullet\ \ $ for $t\in\Bbb G$, define $Q_t:X\x\Bbb G\to X\x\Bbb
G$ by $Q_t(x,g):=(x,gt^{-1})$, then $Q_t\circ T_f=T_f\circ Q_t$.
 \subheading{The associated
invariant action }

The  {\it invariant factor} of $(X\x\Bbb G,\B(X\x\Bbb G),\Bbb
P,T_f)$ is a standard probability space $(\Om,\Cal F,P)=(\Om_{\Cal
I},\Cal F_{\Cal I},P_{\Cal I})$ equipped with a measurable map
$\pi:X\x\Bbb G\to\Om$  such that $\Bbb P\circ\pi^{-1}=P,\ \pi\circ
T_f=\pi$ and $\pi^{-1}\Cal F=\Cal I (T_f)$.

Since  $Q_t\Cal I(T_f)=\Cal I(T_f)$ (because $T_f\circ
Q_t=Q_t\circ T_f$), \f$\bullet\ \ $ $\exists$ a $P$-non-singular
$\Bbb G$-action $\goth p:\Om\to\Om$ so that $\pi\circ Q=\goth
p\circ\pi$. \proclaim{Proposition 1.1}\ \ $(\Om,\Cal F,P,\goth p)$
is ergodic iff $(X,\B,m,T)$ is ergodic.\endproclaim\demo{Proof} If
$(X,\B,m,T)$ is  ergodic, then so is  $(X\x\Bbb G,\B(X\x\Bbb
G),\Bbb P,\<T_f,Q\>)$ where $\<T_f,Q\>$ denotes the $\Bbb
Z_+\x\Bbb G$ action defined by $(n,t)\mapsto T_f^n\circ
Q_t\in\text{\rm Aut}\,(X\x\Bbb G)$. Any $\goth p$-invariant,
measurable function on $\Om$ lifts by $\pi$ to a
$\<T_f,Q\>$-invariant, measurable function on $X\x\Bbb G$, which
is  $\Bbb P$-a.e. constant.
\par Conversely, any $T$-invariant function on $X$ lifts to a $T_f$-invariant
 function on $X\x\Bbb G$ which is also $Q$-invariant and thus the lift of a $\goth p$-invariant, measurable function on $\Om$. If $(\Om,\Cal F,P,g)$ is ergodic this function is constant (a.e.).\hfill\qed\enddemo
\par The $\nsG$ $(\Om,\Cal F,P,\goth p)$ is called the  {\it
invariant-} or {\it Poisson $\Bbb G$-action associated to $(T,f)$}
and denoted $\goth p=\goth p(T,f)$.

 This action is related to   the {\it  Mackey range} of a cocycle (see \cite{Zi2} and \S3),
and the {\it  Poisson boundary} of a random walk (see \S4).\ \
 \subheading{The
associated tail action }

The  {\it tail factor} of $(X\x\Bbb G,\B(X\x\Bbb G),\Bbb P,T_f)$
is a standard probability space $(\Om,\Cal F,P)=(\Om_{\Cal T},\Cal
F_{\Cal T},P_{\Cal T})$ equipped with a measurable map
$\pi:X\x\Bbb G\to\Om$  such that $\Bbb P\circ\pi^{-1}=P,\
\pi^{-1}\Cal F=\Cal T(T_f)$.

 Since $Q_t\Cal T(T_f)=\Cal T(T_f)$ (because $T_f\circ
Q_t=Q_t\circ T_f$), \f$\bullet\ \ $ $\exists$ a $P$-non-singular
$\Bbb G$-action $\tau:\Om\to\Om$ so that $\pi\circ
Q=\tau\circ\pi$. \proclaim{Proposition 1.2}\ \ $(\Om,\Cal
F,P,\tau)$ is ergodic iff $(X,\B,m,T)$ is
exact.\endproclaim\demo{Proof} Suppose first that $(X,\B,m,T)$ is
exact. Any $\tau$-invariant, measurable function
 $F:\Om\to\Bbb R$ lifts by $\pi$ to a $Q$-invariant, $\Cal T(T_f)$-measurable function
  $\o F:X\x\Bbb G\to\Bbb R$. In
particular $\exists\ !\ \o F_n:X\x\Bbb G\to\Bbb R\ \ (n\ge 0)$,
measurable so that $\o F=\o F_n\circ T_f^n$. It follows from
$Q_t\circ T_f=T_f\circ Q_t$ that each $\o F_n$ is $Q$-invariant,
whence $\exists\  \hat F_n:X\to\Bbb R$ measurable so that $\o
F_n(x,y)=\hat F_n(x)$ for $\Bbb P$-a.e. $(x,y)\in X\x\Bbb G$. Thus
$\o F=\o F_0=\hat F_0$  is $\Cal T(T)$-measurable, whence constant
$\Bbb P$-a.e. by exactness of $T$.
\par Conversely, any $\Cal T(T)$-measurable function on $X$ lifts to a
$\Cal T(T_f)$-measurable
 function on $X\x\Bbb G$ which is also $Q$-invariant and thus the lift of a $\tau$-invariant,
 measurable function on $\Om$. If $(\Om,\Cal F,P,\tau)$ is ergodic this function is constant (a.e.).\hfill\qed\enddemo
\par The $\nsG$ $(\Om,\Cal F,P,\tau)$ is called the  {\it associated tail $\Bbb G$-action}
 (of $(T,f)$) and denoted $\tau=\tau (T,f)$. As with the Poisson action, the tail action is
  related to the Mackey range of a cocycle, and also to the {\it tail
 boundary} of a random walk (see \S4).
 \heading{\S2 Conditions for exactness and  a construction of
Zimmer }\endheading\par We begin with a proposition generalising
Zimmer's construction (in \cite{Zi1}) of a $\Bbb G$-valued cocycle
over an $\eppt$ with a prescribed $\ensG$ as Mackey range.
\proclaim{Proposition 2.1}
\par Suppose that  $\Bbb G$ is a LCP group, \f that $(X,\B,m,T)$ is a
non-singular endomorphism, and suppose that $f:X\to\Bbb G$ is
measurable.\par Let $S$ be a non-singular $\Bbb G$-action on a
probability space $(Y,\C,\nu)$ and define $\widetilde T=\widetilde
T_{f,S}:X\x Y\to X\x Y$  by
$$\widetilde T(x,y):=  (Tx,S_{f(x)}y);$$ then
$$\goth p(\widetilde T,f)\cong \goth p(T,f)\x S,\ \&\ \tau (\widetilde T,f)\cong \tau (T,f)\x S.$$
\endproclaim \demo{Proof}
\par Define $\pi:X\x Y\x\Bbb G\to X\x\Bbb G\x Y$ by
$\pi(x,y,g):=(x,g,S_{g^{-1}}y).$ Evidently $\pi$ is a bimeasurable
bijection.\par Fix $p\in \Cal P(\Bbb G),\ p\sim m_{\Bbb G}$. A
calculation shows that
$$(m\x\nu\x p)\circ \pi^{-1}\sim m\x p\x\nu,\tag1$$ indeed
$$\frac{d(m\x\nu\x p)\circ \pi^{-1}}{d(m\x p\x\nu)}(x,g,y)=\frac{d\nu\circ S_g}{d\nu}(y)=:D(g,y).$$
 \par Next, we
claim that
$$\pi\circ\widetilde T_f=(T_f\x\text{\rm Id}|_Y)\circ\pi.\tag2$$ To see
this,
$$\align \pi\circ\widetilde T_f(x,y,g) &=\pi(Tx,S_{f(x)}y,f(x)g)=
(Tx,f(x)g,S_{(f(x)g)^{-1}}S_{f(x)}y)\\ &
=(Tx,f(x)g,S_{g^{-1}}y)=T_f\x\text{\rm
Id}\circ\pi(x,y,g).\endalign$$\par It follows from (2) that
$$\Cal I(\widetilde T_f)=\pi^{-1}(\Cal I(T_f\x\text{\rm Id})),\ \ \Cal T(\widetilde T_f)=\pi^{-1}(\Cal T(T_f\x\text{\rm Id})).$$
 Now, in general
$$\Cal I(T_f\x\text{\rm Id})=\Cal I(T_f)\otimes\B(Y),\ \ \Cal T(T_f\x\text{\rm Id})=\Cal T(T_f)\otimes\B(Y)\ \mod\ m\x p\x\nu$$
and so
$$\Cal I(\widetilde T_f)=\pi^{-1}(\Cal I(T_f)\otimes\B(Y)),\ \ \Cal T(\widetilde T_f)=\pi^{-1}(\Cal T(T_f)\otimes\B(Y))\ \mod\ m\x p\x\nu.$$
 \par Now let $(\Om_i,\Cal F_i,P_i)$ and let $(\widetilde{\Om}_i,\widetilde{\Cal
F}_i, \widetilde P_i)\ \ (i=\goth p,\ \tau)$ be the invariant or
tail factors of $(X\x\Bbb G,m\x p,T_f)$ and $(X\x Y\x\Bbb
G,m\x\nu\x p,\widetilde{T}_f)$ respectively, according to the
value of $i=\goth p,\ \tau$. By (1) and (2), $\pi$ induces a
measure space isomorphism of $(\widetilde{\Om}_i,\widetilde P_i)$
with $(\Om_i\x Y,P_i\x\nu)$.
\par Denoting the associated  $\Bbb G$-actions by  $\widetilde Q_t(x,y,g):=(x,y,gt^{-1})$ and
$Q_t(x,g):=(x,gt^{-1})$, we note that $$\pi\circ \widetilde
Q_{t}=(Q_{t}\x S_t)\circ\pi.$$
 The
proposition now  follows from this. \hfill\qed\enddemo \proclaim{
Corollary 2.2}\par 1) $\widetilde T$ is ergodic iff $\goth
p(T,f)\x S$ is ergodic.\par 2) $\widetilde T$ is exact iff
$\tau(T,f)\x S$ is ergodic.\par 3) If both $T_f$ and $S$ are
ergodic, then $\widetilde T$ is ergodic and $\goth p(\widetilde
T,f)\cong S.$\par 4) If  $T_f$ is exact and $S$ is ergodic, then
$\widetilde T$ is exact and $\tau(\widetilde T,f)\cong S.$
\endproclaim\demo{Proof} Parts 1) and 2) follow from propositions 1.1 and 1.2. Parts
3) and 4) follow from these and form essentially Zimmer's
construction. \hfill\qed\enddemo
\bigskip\heading{ \S3 Locally invertible endomorphisms}\endheading
\par In this section, we obtain additional results for a non-singular, exact endomorphism
  $(X,\B,m,T)$  of a standard measure space which is {\it
locally invertible} in the sense that $\exists$ an at most
countable partition $\a\subset\B$ so that $T:a\to Ta$ is
invertible, non-singular $\forall\ a\in\a$.  Under the assumption
of local invertibility, the associated actions of \S1 are {\it
Mackey ranges of cocycles} (as in \cite{Zi2}). See proposition 3.2
(below).

As  in \cite{S-W}, we call  an $\ensG$ $(X,\B,m,U)$ {\it properly
ergodic} if $m(U_{\Bbb G}(x))=0\ \forall\ x\in X$ and call a
properly   $\ensG$ $S:\Bbb G\to\text{\rm Aut}\,(Y)$ {\it mildly
mixing} if
 $U\x S$ is ergodic for any properly
ergodic $\nsG$  $(X,\B,m,U)$. As shown in \cite{S-W}: \f$\bullet\
\ $ there are no mildly mixing actions of compact groups,
\f$\bullet\ \ $  a  mildly mixing action of non-compact LCP group
has an equivalent, invariant probability. Moreover, \f$\bullet\ \
$
 a $\ppG$ ($\Bbb G$ a non-compact LCP group)
$(Y,\C,\nu,S)$ is mildly mixing iff
$$f\in L^2(\nu),\ g_n\in\Bbb G,\ g_n\overset{\Bbb G}\to\longrightarrow\infty,\ f\circ
S_{g_n}\overset{L^2(\nu)}\to\longrightarrow f\ \Rightarrow\ f\ \
\text{\rm is constant}.$$  We prove

\proclaim{Theorem 3.1}\par Suppose that $(X,\B,m,T)$ is a
non-singular, locally invertible, exact endomorphism of a standard
measure space, that $\Bbb G$ is a LCP, non-compact, Abelian group
and that $f:X\to\Bbb G$ is measurable. \par Either $\widetilde
T_{f,S}$ is exact for every mildly mixing $\ppG$ $S:\Bbb G\to
\text{\rm Aut}\,(Y)$, or $\exists$ a compact subgroup $\Bbb
K\le\Bbb G,\ t\in\Bbb G$ and $\o f:X\to\Bbb K,\ g:X\to\Bbb G$
measurable so that $$f=g-g\circ T+t+\o f.$$

 \endproclaim Note that the invertible version of this generalizes corollary 6 of
\cite{Ru}. The rest of this section is the proof of theorem 3.1.
 \subheading{Tail relations}\par Let
$(X,\B,m,T)$ be a non-singular, locally invertible endomorphism of
a standard probability space. Consider the
 {\it tail relations}
 $$\goth T(T):=\{(x,y)\in X\x X:\ \exists\ k\ge 0,\
T^kx=T^k y\};$$
$$\goth G(T):=\{(x,y)\in X_0\x X_0:\ \exists\ k,\ \ell\ge 0,\ T^kx=T^\ell y\}$$
where $X_0:=\{x\in X:\ T^{n+k}x\ne T^kx\ \forall\ n,\ k\ge 1\}$.
We assume that $m(X\setminus X_0)=0$ (which is the case if $T$ is
ergodic and $m$ is non-atomic) and so $\goth T(T)\subset \goth
G(T)\ \mod m$. Both $\goth T(T)$ and $\goth G(T)$ are {\it
standard, countable, $m$-non-singular equivalence relations} in
sense of \cite{F-M} whose invariant sets are given by
$$\Cal I(\goth
G(T))=\Cal I(T),\ \ \ \Cal I(\goth T(T))=\Cal T(T)$$ respectively.

 Given a LCP group $\Bbb G$ and $f:X\to\Bbb G$ measurable, define
$f_n:X\to\Bbb G\ \ (n\ge 1)$ by
$$f_n(x):=f(T^{n-1}x)f(T^{n-2}x)\dots
f(Tx)f(x)$$ and define $\Psi_f:\goth G(T)\to\Bbb G$ by

$$\Psi_f(x,x'):=f_\ell(x')^{-1}f_k(x)\ \ \ \text{\rm for $k,\ \ell\ge 0$\ \ such that\ \ } T^kx=T^\ell x'$$
(this does not depend on the $k,\ \ell\ge 0$ such that
$T^kx=T^\ell x'$ for $x,x'\in X_0$).

It follows that  $\Psi_f:\goth G(T)\to\Bbb G$ is a (left) {\it
$\goth G(T)$-orbit cocycle} in the sense that
$$\Psi_f(y,z)\Psi_f(x,y)=\Psi_f(x,z)\ \forall\ (x,y),\
(y,z)\in\goth G(T).$$ Note that since $\goth T(T)\subset \goth
G(T)\ \mod m$, the restriction $\Psi_f:\goth T(T)\to\Bbb G$ is a
(left) $\goth T(T)$-orbit cocycle. \subheading{Mackey ranges of
cocycles} Let $\Cal R$ be a countable, standard, non-singular
equivalence relation
 on the standard measure  space $(X,\B,m)$ and let
$\Psi:\Cal R\to\Bbb G$ be a left $\Cal R $-orbit cocycle. It
follows from  theorem 1 in \cite{F-M}, there is a countable group
$\G$ and a $\nsg$ $(X,\B,m,V)$
 so that
 $$\r{=}\Cal R_V:=\{(x,V_\g x):\ \g\in \Gamma ,\ x\in X\}.$$
 Let $f(\g,x):=\Psi(x,V_\g x)\ \ (f=f_{\Psi,V}:\G\x X\to\Bbb G)$ be the associated
  {\it left $V$-cocycle} (satisfying $f(\g\g',z)=f(\g,V_{\g'} z)f(\g',z)$).

The {\it Mackey range}
 $\goth R(V,f)$ (see \cite{Zi2}) is analogous to the invariant action of \S1.
 It is the non-singular  $\Bbb
G$-action of $Q\ \ (Q_t(x,y):=(x,yt^{-1}))$ on the invariant
factor of
 $V_f:X\x\Bbb G\to X\x\Bbb G\ \ ((V_f)_\g(x,g):=(V_\g
x,f(\g,x)g))$.

\bigskip As before, $\goth R(V,f)$ is ergodic iff $(X,\B,m,V)$ is ergodic.

It can be shown that $\goth R(V,f_{\Psi,V})$ does not depend on
$V$ such that $\r{=}\Cal R_V$ and we define the  {\it Mackey range
of $\Psi$ over $\r$} as
 $\goth R(\r,\Psi):=\goth R(V,f_{\Psi,V}).$

\proclaim{Proposition 3.2} $$\goth p(T,f)=\goth R(\goth
G(T),\Psi_f),\ \ \tau(T,f)=\goth R(\goth
T(T),\Psi_f).$$\endproclaim

\bigskip We also have the following version of proposition 2.1
(whose proof is similar) :\proclaim{Proposition 3.3}
\par Suppose that $\G$ is a countable group, that $\Bbb G$ is a LCP group,
\f that $(X,\B,m,V)$ is an $\ensg$, and suppose that $f:\G\x
X\to\Bbb G$ is a measurable cocycle.\par Let $S$ be a non-singular
$\Bbb G$-action on a probability space $(Y,\C,\nu)$ and define
$\widetilde V:\G\to \text{\rm Aut}(X\x Y)$ by $\widetilde
V_\g(x,y):=(V_\g x,S_{f(\g,x)}(y))$; then
$$\goth R(\widetilde V,f)\cong \goth R(V,f)\x S.$$
\endproclaim
\subheading{Compact reducibility}
\par Let $\G$ be a countable
group, $\Bbb G$ be a LCP group and let\f $(X,\B,m,V)$ be an
$\ensg$.

 We call a measurable $V$-cocycle $F:\G\x X\to\Bbb G$ {\it compactly reducible}
 if $\exists\ \Bbb
K\leq \Bbb G$ compact,
 a measurable cocycle  $f:\G\x X\to\Bbb K$ and $h:X\to\Bbb G$ measurable so that
$F(\g,x)=h(V_\g x)^{-1}f(\g,x)h(x)$. \subheading{Regularity and
range of an orbit cocycle} Let $\Bbb G$ be a LCP group and let
$\Cal R\in\B(X\x X)$ be a standard, countable, non-singular
equivalence relation. We call the left $\Cal R $-orbit cocycle
$\Psi:\Cal R\to\Bbb G$ {\it compactly reducible} if the associated
$f=f_{\Psi,V}$ has this property for some (and hence every) $\nsg$
$(X,\B,m,V)$ with
  $\r{=}\Cal R_V$.

  \

For the rest of the section, we assume that $\Bbb G$ is Abelian.
\proclaim{Proposition 3.4} If the measurable $V$-cocycle $F:\G\x
X\to\Bbb G$ is not compactly reducible, then $\goth R(V,F)\x S$ is
ergodic for any mildly mixing \f$\ppG$ $S:\Bbb G\to\text{\rm
Aut}\,(Y)$.
\endproclaim
\demo{Proof} Suppose that the conclusion fails, then $\goth
R(V,F)$
 is not properly ergodic. It follows
from \cite{Zi2}, \ proposition 4.2.24 that    $\exists\ \Bbb H\leq
\Bbb G$ a closed subgroup which is non-compact by assumption,
 a measurable cocycle  $f:\G\x X\to\Bbb H$ and $h:X\to\Bbb G$ measurable so that
$F(\g,x)=f(\g,x)+h(x)-h(V_\g x)$ and such that $(X\x\Bbb
H,\B(X\x\Bbb H),m\x m_{\Bbb H},V_{f})$ is ergodic. In this case,
$\goth R(V,F)$ is the action of $\Bbb G$ on $\Bbb G/\Bbb H$. \par
Let $S:\Bbb G\to\text{\rm Aut}\,(Y)$ be a mildly mixing $\ppG$,
then since $\Bbb H$ is not compact, $S|_{\Bbb H}$ is mildly
mixing, whence ergodic. Also $\goth R(V,F)|_{\Bbb H}=\text{Id.}$
 To see that $\goth R(V,F)\x S$ is ergodic, let $F$ be bounded,
measurable and $\goth R(V,F)\x S$-invariant. For $h\in\Bbb H$,
$F\circ (\goth R(V,F)\x S)_h(\om,y)=F(\om,S_hy)$. By ergodicity of
 $S|_{\Bbb H}$, $\exists\ \o F(\om)$ so that a.e. $F(\om,y)=\o F(\om)$.
 The function $\o F$ is $\goth R(V,F)$-invariant, whence
constant. \hfill\qed\enddemo

 By proposition 3.4,
\proclaim{Proposition 3.5}
\par If $\Psi$ is not compactly
reducible, then $\goth R(\r,\Psi)\x S$ is ergodic for any mildly
mixing \f$\ppG$ $S:\Bbb G\to\text{\rm Aut}\,(Y)$. \endproclaim
  We now have by propositions 3.2 and 3.5 that
\proclaim{ Proposition 3.6} \par  For a locally invertible, exact
 endomorphism $(X,\B,m,T)$ and $f:X\to\Bbb G$ measurable:
  If $\Psi_f$ is not
$\goth T(T)$-compactly reducible, then $\tau(T,f)\x S$ is ergodic
for any mildly mixing \f$\ppG$ $S:\Bbb G\to\text{\rm Aut}\,(Y)$.
\endproclaim \subheading{Proof of theorem 3.1}
\par The previous propositions show that  if $\Psi_f$ is not $\goth T(T)$-compactly
reducible, then $\widetilde T_{f,S}$ is exact for every mildly
mixing \f$\ppG$ $S:\Bbb G\to \text{\rm Aut}\,(Y)$.

We must show that if $\Psi_f$ is $\goth T(T)$-compactly reducible
then $\exists$ a compact subgroup $\Bbb K\le\Bbb G,\ t\in\Bbb G$
and $\o f:X\to\Bbb K,\ g:X\to\Bbb G$ measurable so that
$$f=g-g\circ T+t+\o f.$$

To see this suppose that  $\Bbb K\le\Bbb G$ is a compact subgroup
and $g:X\to\Bbb G$ is measurable so that the quotient cocycle
(under the map $s\mapsto\widetilde s:=s +\Bbb K\ \ (\Bbb G\to\Bbb
G/\Bbb K)$) is a coboundary, i.e.
$\widetilde\Psi_f(x,x')=\widetilde g(x)-\widetilde g(x')$. Set
$h(x):=\widetilde f(x)-\widetilde g(x)+\widetilde g(Tx)$.

We claim that  $h$ is $\goth T(T)$-invariant, whence constant.

To prove this $\goth T(T)$-invariance, (as in proposition 1.2 of
\cite{ANS}), note that
 $T^nx=T^nx'\ \Rightarrow\ h_n(x)-h_n(x')=\widetilde f_n(x)-\widetilde f_n(x')-\widetilde
g(x)+\widetilde g(x')=0$. Also $T^nx=T^nx'\ \Rightarrow\
T^{n-1}(Tx)=T^{n-1}(Tx')\ \Rightarrow\ h_{n-1}(Tx)=h_{n-1}(Tx')$.
Since $h(x)=h_n(x)-h_{n-1}(Tx)$, we have $h(x)=h(x')$. \hfill\qed

 \heading{\S4 Random walks}\endheading
Let $\Bbb G$ be a LCP group and let $p\in\Cal P(\Bbb G)$.

 The (left) {\it random walk on $\Bbb G$} with
{\it jump probability } $p\in\Cal P(\Bbb G)$ ($\text{\tt RW}(\Bbb
G,p)$ for short) is the stationary, one-sided shift of the Markov
chain on $\Bbb G$ with transition probability $P(g,A):=p(Ag^{-1})\
\ (A\in\B(\Bbb G))$. The random walk $\text{\tt RW}(\Bbb G,p)$ is
said to be {\it adapted} if   $p\in\Cal P(\Bbb G)$ is {\it
globally supported} in the sense that $\o{\<\text{\rm
supp}\,(p)\>}=\Bbb G$.

\smallskip\f The random walk $\text{\tt RW}(\Bbb G,p)$ is isomorphic to the $\mpt$ $$(X\x\Bbb G,\B(X\x\Bbb
G),\mu_p\x m_{\Bbb G},W)$$ where $X=\Bbb G^{\Bbb N}$,
$\mu_p:\B(X)\to [0,1]$ is the product measure $\mu_p:=p^\Bbb N$
defined by $\mu_p([A_1,\dots,A_n]):=\prod_{k=1}^np(A_k)$ where
$[A_1,\dots,A_n]:=\{x\in X:\ x_k\in A_k\ \forall\ 1\le k\le n\}$,
($n\ge 1,\ A_1,\dots,A_n\in\B(\Bbb G)$), $m_{\Bbb G}$ is left Haar
measure on $\Bbb G$ and $W:X\x\Bbb G\to X\x\Bbb G$ is defined by
$W(x,g):=(Tx,x_1g)$ with $T:X\to X$ being the shift
$(Tx)_n:=x_{n+1}$.

The {\it jump process} of the random walk $\rwgp$ is
$$(X,\B(X),\mu_p,T,f)$$ where $(X,\B(X),\mu_p,T)$ is as above and
$f:X\to\Bbb G$ is defined by $f(x):=x_1$.

\subheading{Boundaries}

The {\it tail boundary} of the random walk $\text{\tt RW}(\Bbb
G,p)$ is $\tau(\Bbb G,p):=\tau(T,f)$, and the {\it Poisson
boundary} is $\goth p(\Bbb G,p):=\goth p(T,f)$ where
$(X,\B(X),\mu_p,T,f)$ is the {\it jump process} of the random walk
$\rwgp$.

These definitions are equivalent with those in \cite{K-V,\ K} (see
also \cite{Fu}) . \subheading{Weakly mixing actions}
\par Let $\Bbb G$ be a LCP group. A $\nsG$ $V:\Bbb G\to\text{\rm Aut}\,(X,\B,m)$ is called {\it weakly mixing}
if $V\x S$ is ergodic on $X\x Y$ whenever $S:\Bbb G\to\text{\rm
Aut}\,(Y,\Cal C,\nu)$ is an $\eppG$. \par In case $\Bbb G=\Bbb Z$,
this agrees with the definition given in \cite{ALW}.  More
generally, in case $\Bbb G$ is Abelian, weak mixing of $V$ is
equivalent to the condition
$$f\in L^\infty(X,\B,m),\ \g\in\widehat{\Bbb G},\ f\circ
V_g=\g(g)f\ \text{a.e.}\ \forall\ g\in\Bbb G\ \Rightarrow\ \ f\ \
\text{is a.e. constant}$$ where $\widehat{\Bbb G}$ denotes the
dual group of $\Bbb G$. For a proof of this equivalence when $\Bbb
G=\Bbb Z$, see \S4 of \cite{ALW}. Alternatively, see theorem 2.7.1
in \cite{A} (which can easily be extended to the general Abelian
case). \subheading{Random walks on Abelian groups} In case $\Bbb
G$ is Abelian, the tail and Poisson boundaries are the actions of
$\Bbb G$ on $\Bbb G/\Bbb H_{\tau }$ and $\Bbb G/\Bbb H_{\goth p}$
by translation (respectively), where
$$\Bbb H_{\goth p}:=\o{\<\text{\rm supp}\,p\>},\ \Bbb H_{\tau }:=
\o{\<\text{\rm supp}\,p\ -\  \text{\rm supp}\,p\>}.$$ See
\cite{D-L}. Here, for $F\subset\Bbb G,\ \<F\>$ denotes the minimal
subgroup of $\Bbb G$ containing $F$.

 \proclaim{Theorem 4.1\ \ \ (extension of \cite{Me})}\par
 Let $\Bbb G$ be a LCP, Abelian group, let $p\in\Cal
P(\Bbb G)$ and let
 $(X,\B(X),\mu_p,T,f)$ be
the  jump process of the random walk $\rwgp$. \par 1)\ \ If \ \
$\Bbb G/\Bbb H_{\goth p}\ \ (\Bbb G/\Bbb H_{\tau})$\ \ is compact
and $S:\Bbb G\to\text{\rm Aut}\,(Y)$ is a weakly mixing $\ppG$,
then $\widetilde T_{f,S}$ is ergodic (exact).
\par 2)\ \ If \ \
$\Bbb H_{\goth p}\ \ (\Bbb H_{\tau})$\ \ is non-compact and
$S:\Bbb G\to\text{\rm Aut}\,(Y)$ is a mildly mixing $\ppG$, then
$\widetilde T_{f,S}$ is ergodic (exact).\endproclaim \demo{Proof}
\f 1) It follows from the assumption that $\goth p(T,f)\ \ (
\tau(T,f))$ is an \f$\eppG$, whence
 $S\x \goth p(T,f)\ \ (S\x \tau(T,f))$ is ergodic whenever $S$ is weakly mixing.
 \f 2) The proof is  as  the proof of proposition 3.4 (above).\hfill\qed\enddemo

 \proclaim{Proposition 4.2\ \ \ (converse of \cite{Me})}\par Let
$p\in\Cal P(\Bbb Z)$ and let
 $(X,\B(X),\mu_p,T,f)$ be
the  jump process of the random walk $\text{\rm RW}\,(\Bbb Z,p)$.
Suppose
$$\Bbb H:=\o{\<\text{\rm supp}\,(p)-\text{\rm supp}\,(p)\>}=d\Bbb
Z$$ with $d\ne\{0\}$ and let $S$ be an $\eppt$ of $Y$, then
 $\widetilde T_{f,S}$ is exact iff $S^d$ is ergodic.\endproclaim \demo{Proof}
Here,  $\tau(T,f)$ is the cyclic permutation of $d$ points. The
ergodicity of $S^d$    characterizes the ergodicity of
$S\x\tau(\goth T(T),\psi_f)$. \hfill\qed\enddemo

\subheading{Remark} \ \ In case $\Bbb G$ is LCP, Abelian,
non-compact and\ \ $\Bbb H_i$\ \ is compact ($i=\goth p$ or
$\tau$), then a Gaussian action of $\Bbb G/\Bbb H_i$ with Haar
spectral measure type  is mildly mixing. Lifting this action, we
obtain a mildly mixing $\ppG$
 $S:\Bbb G\to\text{\rm Aut}\,(Y)$  with $S|_{\Bbb H_i}\equiv\text{Id}.$
  Here,
$\widetilde T_{f,S}$ is not ergodic or exact (according to whether
$i=\goth p$ or $\tau$). This is because  the action $V\x S$ on
$\Bbb G/\Bbb H_i\x Y$ is not ergodic (where $V_t(g+\Bbb
H_i):=t+g+\Bbb H_i$).

Indeed, if $A\in\B(Y)$, then $S_kA=A\ \mod\ \nu\ \forall\ k\in\Bbb
H_i$. It follows that $B=B(A):=\bigcup_{g\in\Bbb G}g\Bbb H_i\x gA$
is Lebesgue measurable, $V\x S$-invariant, $m_{\Bbb G/\Bbb
H_i}\x\nu(B(A))>0$ for $\nu(A)>0$ and  $B(A)\cap
B(A^c)=\emptyset$.\subheading{Weak mixing of Poisson boundary}

We show that the Poisson action $\goth p(\Bbb G,p)$ is weakly
mixing when  $\text{\rm RW}\,(\Bbb G,p)$ is adapted.

Let $p\in\Cal P(\Bbb G)$ is globally supported and let
$(X,\B(X),\mu_p,T,f)$ be the jump process of the random walk
$\text{\rm RW}\,(\Bbb G,p)$. \proclaim{Proposition 4.3}\par
 Suppose that $S:\Bbb G\to\text{\rm Aut}\,(Y)$ is a
 \f $\ppG$, then $\widetilde T_{f,S}$ is ergodic iff $S$ is ergodic.
\endproclaim \demo{Proof of $S$ ergodic $\Rightarrow\
 \widetilde T_{f,S}$  ergodic {\rm as in \cite{Mo}, (see also \cite{ADSZ}, proposition 1)}}

 Suppose that $h:X\x Y\to\Bbb R$ is
 bounded, measurable and $\widetilde T_{f,S}$-invariant, then $P_{\widetilde T_{f,S}}h=h$
 where $P_{\widetilde T_{f,S}}:L^1(\mu_p\x\nu)\to L^1(\mu_p\x\nu)$
 is the predual of $$F\mapsto F\circ \widetilde T_{f,S}\ \ (L^\infty(\mu_p\x\nu)\to
 L^\infty(\mu_p\x\nu))$$ and given by
$$P_{{\widetilde T}_{f,S}}^nF(x,y)=P_T^n(F(\cdot,S_{\a_n(\cdot)}^{-1}(y))(x)$$
where $\a_n(x):=x_n\dots x_1$ and $P_T:L^1(\mu_p)\to L^1(\mu_p)$
 is the predual of $F\mapsto F\circ T\ \ \ \ (L^\infty(\mu_p)\to
 L^\infty(\mu_p))$.

 We claim that $h$ is $X\x\Cal C$-measurable.

To see this,  note that
 $\exists\ h_n,\ \s(x_1,\dots,x_n)\x\Cal
C$-measurable so that $\|h-h_n\|_1\to 0$. By independence of
$x_1,x_2,\dots$,
$$\align &
P_{\widetilde T_{f,S}^n}h_n(x,y)=P_T^n(h_n(\cdot,S_{\a_n(\cdot)}^{-1}(y))(x)\\
& =E(h_n(\cdot,S_{\a_n(\cdot)}^{-1}(y))= E(P_{\widetilde
T_{f,S}^n}h_n|X\x\Cal C).\endalign$$ Thus

$$\align  &\|h-E(h|X\x\Cal C)\|_1\\ &\le \|P_{\widetilde
T_{f,S}^n}h-P_{\widetilde T_{f,S}^n}h_n\|_1 +\|E(P_{\widetilde
T_{f,S}^n}h_n|X\x\Cal C)-E(P_{\widetilde T_{f,S}^n}h|X\x\Cal
C)\|_1\\ & \leq 2\|h-h_n\|_1\to 0.
\endalign$$

 Thus $h=G$ where $G:Y\to\Bbb R$ and $G\circ S_g=G\ \nu$-a.e. $\forall\ g\in\o{\<\text{\rm
supp}\,(p)\>}=\Bbb G$ and $G$ (whence $h$) is a.e. constant by
ergodicity of $S$.   \hfill\qed\enddemo \proclaim{Proposition
4.4}\par $\goth p(T,f)$ is weakly mixing in the sense that $\goth
p(T,f)\x S$ whenever $S$ is an\f $\eppG$.
\endproclaim\demo{Proof} Let $S:\Bbb G\to\text{\rm Aut}\,(Y)$ be
an  $\eppG$. By proposition 4.3,  $\widetilde T_{f,S}$ is ergodic,
whence by proposition  2.1,
 $\goth p(T,f)\x S$ is ergodic. \hfill\qed\enddemo
\subheading{Remark}

A {\it fibred system} $(X,\B,m,T,\a)$ is a non-singular
endomorphism $(X,\B,m,T)$ which is locally invertible with respect
to the  at most countable partition $\a\subset\B$, which also
satisfies $\s(\{T^{-n}a:\ n\ge 0,\ a\in\a\})\overset{m}\to{=}\B$.
Propositions 4.3 and 4.4 remain true whenever $(X,\B,m,T,\a)$ is a
  probability preserving fibred system, which is {\it weak
quasi--Markov, almost onto} in the sense of \cite{ADSZ} and
$f:X\to\Bbb G$ is $\a$-measurable with  $m\circ f^{-1}$ is
globally supported on $\Bbb G$.
 \subheading{Aperiodic random walks}
\par Let $\Bbb G$ be a LCP group and suppose that $p\in\Cal
P(\Bbb G)$.

A random walk $\rwgp$ is called {\it steady} (see \cite{K}) if
$\goth p(\Bbb G,p)=\tau(\Bbb G,p)$.

\par Let $\Bbb G$ be a countable group and suppose that $p\in\Cal
P(\Bbb G)$.
 The
 random walk $\text{\tt RW}(\Bbb G,p)$  is called {\it aperiodic} if the corresponding Markov
chain is aperiodic. Equivalent conditions for this are \f$\bullet\
\ \ $ $p^{n*}_e>0 \ \forall\ n$ large; \f$\bullet\ \ \ $ $\exists\
n\ge 1$ so that $\text{\rm supp}\,(p^{n*})\cap \text{\rm
supp}\,(p^{(n+1)*})\ne\emptyset$.

 An aperiodic random walk on a countable group is steady.
This can be gleaned from \cite{Fo} and th\'eor\`eme 3 in \cite{D}
(see also proposition 4.5 in \cite{K}). \proclaim{Theorem 4.5}\par
Let $\Bbb G$ be a LCP group and suppose that $p\in\Cal P(\Bbb G)$
is  globally supported and that $\rwgp$ is  steady. Let
$(X,\B(X),\mu_p,T,f)$ be the jump process of the random walk
$\text{\rm RW}\,(\Bbb G,p)$ and let $S:\Bbb G\to\text{\rm
Aut}\,(Y)$ be a
 $\ppG$, then $\widetilde T_{f,S}$ is ergodic iff $S$ is ergodic and in this case
 $\widetilde T_{f,S}$ is exact.\endproclaim \demo{Proof {\rm of $S$ ergodic $\Rightarrow\
 \widetilde T_{f,S}$  exact}}   By proposition 4.4,
 $\goth p(T,f)\x S$ is ergodic, whence by steadiness,
  $\tau(T,f)\x S$ is ergodic.
Thus, $\widetilde T_{f,S}$ is  exact by corollary 2.2, 2).
\hfill\qed\enddemo\heading References\endheading  \enddocument